# Uniquely universal sets in $\mathbb{R} \times \omega$ and $[0,1] \times \omega$

Alicja Krzeszowiec

**Abstract**

Let $X$ and $Y$ be topological spaces. We say that $X \times Y$ satisfies the Uniquely Universal property (UU) iff there exists an open set $U \subseteq X \times Y$ such that for every open set $W \subseteq Y$ there is a unique cross section of $U$ with $U(x) = \{y \in Y : (x, y) \in U\} = W$. Arnold W. Miller in his paper [1] posed the following two questions:

**Problem 1** *Does $[0,1] \times \omega$ have UU?*

**Problem 2** *Does $\mathbb{R} \times \omega$ have UU?*

In this paper we present two constructions which give positive answers to both problems.

Let us first introduce the following notation.

Let $X, Y$ be any spaces, let $x \in X$, $y \in Y$ and $A \subseteq X \times Y$. Then let us denote by $A(x) = \{a \in Y : (x, a) \in A\}$ the horizontal section of set $A$, and by $(y)A = \{b \in X : (b, y) \in A\}$ the vertical section of set $A$.

Let $x \in 2^n$, $x = (x(0), x(1), x(2), ..., x(n-1))$, and let $i \in \{0, 1\}$. By $x \frown i$ we denote the sequence $(x(0), x(1), ..., x(n-1), i)$. Let $k < n$, by $x_{|k}$ we denote the sequence consisting of first $k$ elements of sequence $x$, i.e. $x_{|k} = (x(0), x(1), ..., x(k-1))$.

**Theorem 3** *There is UU set for $[0,1] \times \omega$.*

**Proof.** We will find uniquely universal closed set. The construction is as follows: Let us divide the interval $[0,1]$ into two intervals of equal length $A_0 = \left[0, \frac{1}{2}\right)$, $A_1 = \left[\frac{1}{2}, 1\right]$. In the next step we divide each $A_0$ and $A_1$ into two intervals of equal length, $A_{01} = \left[0, \frac{1}{4}\right]$, $A_{00} = \left(\frac{1}{4}, \frac{1}{2}\right)$, $A_{10} = \left[\frac{1}{2}, \frac{3}{4}\right)$, $A_{11} = \left[\frac{3}{4}, 1\right]$. And we proceed, i.e. in the $n$th step we divide intervals into two equal parts, and looking from the bottom: the first one, $A_{x_0 \frown 1} = \left[0, \frac{1}{2^{n+1}}\right]$, $A_{x_0 \frown 0} = \left(\frac{1}{2^{n+1}}, \frac{2}{2^{n+1}}\right)$, $A_{x_1 \frown 0} = \left[\frac{2}{2^{n+1}}, \frac{3}{2^{n+1}}\right)$, $A_{x_1 \frown 1} = \left[\frac{3}{2^{n+1}}, \frac{4}{2^{n+1}}\right]$, $A_{x_2 \frown 1} =$ $\left[\frac{4}{2^{n+1}}, \frac{5}{2^{n+1}}\right]$, $A_{x_2 \frown 0} = \left(\frac{5}{2^{n+1}}, \frac{6}{2^{n+1}}\right)$, $A_{x_3 \frown 0} = \left[\frac{6}{2^{n+1}}, \frac{7}{2^{n+1}}\right)$, $A_{x_3 \frown 1} = \left[\frac{7}{2^{n+1}}, \frac{8}{2^{n+1}}\right]$ and so on, where $x_0 = \left(0, \underbrace{1, 1, ..., 1}_{\text{n-4 times}}, 1, 1\right)$, $x_1 = \left(0, \underbrace{1, 1, ..., 1}_{\text{n-4 times}}, 1, 0\right)$, $x_2 = \left(0, \underbrace{1, 1, ..., 1}_{\text{n-4 times}}, 0, 0\right)$. Notice that in each step intervals which are numerated by any sequence finishing with 1 are closed. Below the scheme of the construction is given.

$$
\begin{array}{ccc}
 & & A_{111} \\
 & A_{11} & A_{110} \\
A_1 & & A_{100} \\
 & A_{10} & A_{101} \\
[0,1] & & \\
 & A_{00} & A_{001} \\
A_0 & & A_{000} \\
 & A_{01} & A_{010} \\
 & & A_{011}
\end{array}
$$

Let
$$W = \bigcup_{n \in \omega} \bigcup_{x \in 2^n} A_{x \frown 1} \times \{n\}.$$

As we noticed before, for each $x \in 2^n$ set $A_{x \frown 1}$ is closed, hence $W$ is closed. We have the following properties.

**Property 1.** *Let $\overline{x} \in 2^\omega$ be such that for each $k < \omega$ there is $n \geq k$ such that $\overline{x}(n) = 0$. Then $\bigcap_{n>0} A_{\overline{x}_{|n}} = \{a_{\overline{x}}\}$.*

Indeed, the diameters of $A_{\overline{x}_{|n}}$'s (for any $\overline{x} \in 2^\omega$) tend to zero, hence the intersection has at most one point. If for each $k \in \omega$ there is $n_k > k$ such that $\overline{x}(n_k) = 0$, then the distance between endpoints of $A_{\overline{x}_{|n_k}}$ and the endpoints of $A_{\overline{x}_{|n_k-2}}$ is $\frac{1}{2^n}$. Thus the endpoints of $A_{\overline{x}_{|n_k}}$ belong to the interior of $A_{\overline{x}_{|n_k-2}}$. Hence $\bigcap_{n>0} A_{\overline{x}_{|n}} \neq \emptyset$.

Notice that if $\overline{x}(n) = 1$ for all $n > k$ for some $k$, then it may happen that $\bigcap_{n>0} A_{\overline{x}_{|n}} = \emptyset$ (for example if $\overline{x} = (0, 0, 1, 1, 1, ...)$).

**Property 2.** *For distinct $\overline{x}, \overline{y} \in 2^\omega$ we have $\bigcap_{n>0} A_{\overline{x}_{|n}} \neq \bigcap_{n>0} A_{\overline{y}_{|n}}$ (if nonempty).*

Indeed, let $n$ be the smallest number such that $\overline{x}_{|n} \neq \overline{y}_{|n}$. Without loss of generality we may assume that $\overline{x}_{|n} = \overline{x}_{|n-1} \frown 1$ and $\overline{y}_{|n} = \overline{x}_{|n-1} \frown 0$. But then $A_{\overline{x}_{|n-1}\frown 0} \cap A_{\overline{x}_{|n-1}\frown 1} = \emptyset$. Hence $\bigcap_{n>0} A_{\overline{x}_{|n}} \neq \bigcap_{n>0} A_{\overline{y}_{|n}}$.

**Property 3.** *For each $n > 0$*

$$(n)\, W = \left[0, \frac{1}{2^{n+1}}\right] \cup \left[\frac{2^{n+1}-1}{2^{n+1}}, 1\right] \cup \bigcup_{k=1}^{2^{n-1}-1} \left[\frac{4k-1}{2^{n+1}}, \frac{4k+1}{2^{n+1}}\right],$$

*and* $(0)\, W = \left[\frac{1}{2}, 1\right]$.

Let $A \subseteq \omega$ be any subset of $\omega$. We can identify this subset with sequence $x^A \in 2^\omega$. Let us assume that $x^A$ is such that for each $k < \omega$ there is $n \geqslant k$ such that $x^A(n) = 0$. From Property 1 it follows that $\bigcap_{n \in \omega} A_{x^A_{|n}} = \{a_{x^A}\}$. Hence the section $W(a_{x^A}) = \{n \in \omega : (a_{x^A}, n) \in U\} = A$. Indeed, $(a_{x^A}, n) \in W$ iff $a_{x^A} \in \bigcup_{x \in 2^n} A_{x \frown 1}$. There are at most two sequences $y, x \in 2^n$ such that $a_{x^A} \in A_{x \frown 1}$ and $a_{x^A} \in A_{y \frown 1}$. But then $\bigcap_{k \leqslant n} A_{x_{|k}} \cap \bigcap_{k \leqslant n} A_{y_{|k}} = \emptyset$, hence $a_{x^A}$ belongs only to one of the intersections, say $\bigcap_{k \leqslant n} A_{x_{|k}}$. Then $x = x^A_{|n}$. From the fact that $A_{x \frown 1} \cap A_{x \frown 0} = \emptyset$, it follows that $x \frown 1 = x^A_{|n+1}$. Hence if $(a_{x^A}, n) \in W$, and it follows that $x^A(n) = 1$. The other implication is obvious.

From Property 2 it follows that $x^A$ is the only sequence such that $\bigcap_{n>0} A_{x^A_{|n}} = \{a_{x^A}\}$.

From above it follows that as horizontal sections of set $W$ we obtain in a unique way all sequences which are not constantly equal to 1 from some point on. We also obtain (in a unique way) some sequences which are constantly equal to 1 from some point on, for example, the sequence $(1, 0, 1, 1, 1...)$ is obtained as $W_{\frac{1}{2}}$.

The graph below presents the set $W$.

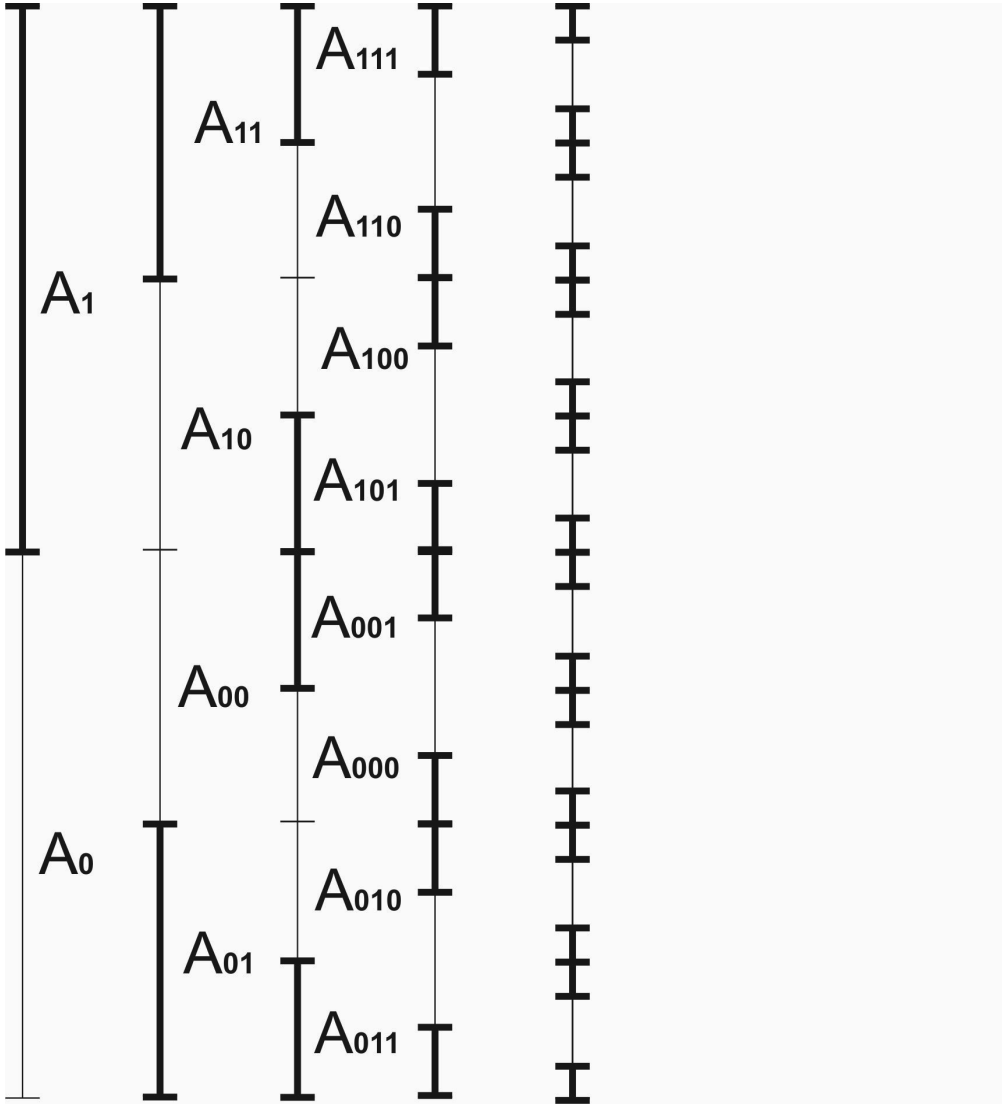

Figure 1:

We need to refine our construction to obtain set $U$ which will be uniquely universal. As it was mentioned above we are missing some of the sequences which are constantly equal to 1 from some point on. Let us notice that when $\bar{x} \in 2^\omega$ such that there is some $k < \omega$ such that $\bar{x}(k) = 0$ and $\bar{x}(n) = 1$ for all $n > k$, then $\bigcap_{n>k+1} A_{\bar{x}|n}$ =common endpoint of $A_{\bar{x}_{|k-1}\frown 0,0}$ and $A_{\bar{x}_{|k-1}\frown 1,0}$ which is the same as the common endpoint of $A_{\overline{x}_{|k+1}}$ and $A_{\overline{x}_{|k-1}\frown 1,0}$. Hence $\bigcap_{n>k+1} A_{\overline{x}_{|n}} \times \{k\} \notin W$. Moreover $\bigcap_{n>k+1} A_{\overline{x}_{|n}}$ is one of the endpoints of interval $A_{\overline{x}_{|k}}$ and $\bigcap_{n>k+1} A_{\overline{x}_{|n}} \in Int\left(A_{\overline{x}_{|t}}\right)$ for all $t < k$. From these properties it follows that if $\overline{x} \in 2^\omega$ is such that there is some $k < \omega$ such that $\overline{x}(k-1) = 1$, $\overline{x}(k) = 0$ and $\overline{x}(n) = 1$ for all $n > k$, then $\bigcap_{n>k+1} A_{\overline{x}_{|n}} = \{a\}$ and the horizontal section $W_a$ generates sequence $\overline{x}$. It also follows that the only sequences which are not generated as horizontal sections of $W$ are sequences $\overline{x} \in 2^\omega$ is such that there is some $k < \omega$ such that $\overline{x}(k-1) = 0$, $\overline{x}(k) = 0$ and $\overline{x}(n) = 1$ for all $n > k$. But there are only countably many such sequences. Let $\{\overline{x_n} : n < \omega\}$ be the list of all such sequences. We will define set $U$ which will be $UU$ by mathematical induction.

Let $k_0 < \omega$ be such that $\overline{x_0}(k_0 - 1) = 0$, $\overline{x_0}(k_0) = 0$ and $\overline{x_0}(n) = 1$ for all $n > k_0$. Let $\overline{y_0}$ be such that $\overline{y_0}(n) = \overline{x_0}(n)$ for all $n \neq k_0 + 2$, and $\overline{y_0}(k_0 + 2) = 0$. Let us first notice that $\overline{y_0} \neq \overline{x_n}$ for all $n < \omega$. Let $\{a_0\} = \bigcap_{n>k_0+3} A_{\overline{y_0}_{|n}} = \bigcap_{n>0} A_{\overline{y_0}_{|n}}$. By adding point $(a_0, k_0 + 2)$ to the set $W$, we obtain sequence $\overline{x_0}$ as horizontal section of $W \cup (a_0, k_0 + 2)$, but we loose sequence $\overline{y_0}$.

Let $U_0 = W \cup (a_0, k_0 + 2)$.

In the next step we add sequence $\overline{y_0}$. Let $\overline{y_1} \in 2^\omega$ be such that $\overline{y_0}(n) = \overline{y_1}(n)$ for all $n \neq k_0+4$ and $\overline{y_1}(k_0 + 4) = 0$. Let $\{a_1\} = \bigcap_{n>k_0+5} A_{\overline{y_1}_{|n}}$. Similarly, by adding point $(a_1, k_0 + 4)$ to set $U_0$ we will obtain sequence $\overline{y_0}$, but we loose sequence $\overline{y_1}$.

Let $U_1 = U_0 \cup (a_1, k_0 + 4)$.

In the next step we add sequence $\overline{x_1}$. Let $k_1 < \omega$ be such that $\overline{x_1}(k_1 - 1) = 0$, $\overline{x_1}(k_1) = 0$ and $\overline{x_1}(n) = 1$ for all $n > k$. Let $\overline{y_2} \in 2^\omega$ be such that $\overline{y_2}(n) = \overline{x_1}(n)$ for all $n \neq k_1+2$, and $\overline{y_2}(k_1 + 2) = 0$. Notice that $\overline{y_2} \neq \overline{y_0}$ and $\overline{y_2} \neq \overline{y_1}$. Indeed, if $\overline{y_2} = \overline{y_0}$, then $\overline{x_1} = \overline{x_0}$ - a contradiction. If $\overline{y_2} = \overline{y_1}$ then this sequence would be of form $(\overline{x_0}(0), \overline{x_0}(1), ..., \overline{x_0}(k_0 - 2), 0, 0, 1, 0, 1, 0, 1, 1, 1...)$, and on the other hand $\overline{y_2} = (\overline{x_1}(0), \overline{x_1}(1), ..., \overline{x_1}(k_1 - 2), 0, 0, 1, 0, 1, 1, ...)$ - a contradiction. Hence $\overline{y_2} \neq \overline{y_0}$ and $\overline{y_2} \neq \overline{y_1}$ and we proceed in the same way as in the first step.

In the next step we add sequence $\overline{x_2}$, and so on. In each of the steps we add one sequence and we loose one sequence, but each sequence is added at some step.

Let $U = \bigcup_{n<\omega} U_n$. Then set $U$ is closed, because the set $W$ was closed and to each of the vertical sections $(n)W$ we add only finitely many points. From the properties of the set $W$ and from the construction it follows that $U$ is uniquely universal for $[0,1] \times \omega$. ∎

**Theorem 4** *There is UU set for $\mathbb{R} \times \omega$.*

**Proof.** We will construct the uniquely universal set in interval $(-1, 1)$. We will divide intervals $(-1, 0)$ and $(0, 1)$ into $\omega$ subintervals. Let $U_{-n}$ be the set $U$ from previous theorem, but constructed in interval $\left(\frac{-1}{n+1}, \frac{-1}{n+2}\right]$ and starting the construction on the $2n + 1$'th coordinate, i.e. $(k) U_{-n} = \emptyset$ for $k \leqslant 2n$ and $(2n+1) U_{-n} = \left[-\frac{2n+3}{(n+1)(n+2)}, -\frac{1}{n+2}\right]$. The only sequence which we will not obtain as horizontal section of set $U_{-n}$ is $\left(\underbrace{0, 0, ...0}_{2n \text{ times}}, 0, 1, 1, 1, 1, ...\right)$, because this sequence can only be obtained as section $\left(\frac{-1}{n+1}\right) U_{-n}$, but $\frac{-1}{n+1} \notin \left(\frac{-1}{n+1}, \frac{-1}{n+2}\right]$. Let $U^-$ be the set constructed in previous theorem, and then reflected with respect to 0 (i.e. $U^- = -U$). Let $U_n$ be the set $U^-$ constructed in interval $\left[\frac{1}{n+2}, \frac{1}{n+1}\right)$ and starting the construction on the $2n + 2$ coordinate (similarly, the only sequence which we will not obtain as horizontal section of $U_n$ is $\left(\underbrace{0, 0, ...0}_{2n+1 \text{ times}}, 0, 1, 1, 1, 1, ...\right)$).

Let

$$A = \left(\bigcup_{n \in \omega} \left(\left(\left(\frac{-1}{n+1}, \frac{-1}{n+2}\right] \times \{2n\}\right) \cup U_{-n}\right)\right)$$
$$\cup \left(\bigcup_{n \in \omega} \left(\left(\left[\frac{1}{n+2}, \frac{1}{n+1}\right) \times \{2n+1\}\right) \cup U_n\right)\right).$$

For $n = 0$ we have $\left(\left(-1, \frac{-1}{2}\right] \times \{0\}\right) \cup U_{-0}$. Horizontal sections of this part of the set $A$ will generate all sequences $x \in 2^\omega$ with $x(0) = 1$, except the sequence $(1, 0, 1, 1, 1, 1...)$. As horizontal sections of part $\left(\left[\frac{1}{2}, 1\right) \times \{1\}\right) \cup U_0^-$ of set $A$ we obtain all sequences $x \in 2^\omega$ such that $x(0) = 0$ and $x(1) = 1$, and the only sequence which we will not obtain is $(0, 1, 0, 1, 1, 1, 1...)$.
For $n = 1$ we get $\left(\left(\frac{-1}{2}, \frac{-1}{3}\right] \times \{2\}\right) \cup U_{-1}$ and $\left(\left[\frac{1}{2}, \frac{1}{3}\right) \times \{3\}\right) \cup U_1$. The first set generates all sequences such that $x(0) = 0, x(1) = 0$ and $x(2) = 1$ (except sequence $(0, 0, 1, 0, 1, 1, 1, ...)$), and so on. Hence as horizontal sections of set $A$ we obtain all sequences except $\left(0, ..., 0, \underset{n'\text{th coordinate}}{1}, 0, 1, 1, 1, 1...\right)$ for all $n \in \omega$. All other sequences are obtained in a unique way.

The set $A$ is closed. Indeed, we will give the proof for the negative part for the set $A$, i.e. that $A \cap (-1, 0]$ is closed, as proof for the positive part is analogical.

Let us first notice that, by Property 3 of set $W$ from the previous theorem, it follows that $\bigcup_{k>0} \left(\frac{-1}{n+1}, \frac{-1}{n+1} + \frac{1}{2^{k+1}(n+1)(n+2)}\right] \times \{2n+1+k\} \subseteq U_{-n}$ and

$$\bigcup_{k=0} \left[\frac{-1}{n+2} - \frac{1}{2^{k+1}(n+1)(n+2)}, \frac{-1}{n+2}\right] \times \{2n+1+k\} \subseteq U_{-n}$$

for each $n \in \omega$. We only have to show that although set $A$ contains half open intervals $\left(\frac{-1}{n+1}, \frac{-1}{n+2}\right] \times \{2n\}$ and $\bigcup_{k>0} \left(\frac{-1}{n+1}, \frac{-1}{n+1} + \frac{1}{2^{k+1}(n+1)(n+2)}\right] \times \{2n+1+k\}$, it is closed.

For $n = 0$ both $\left(\frac{-1}{n+1}, \frac{-1}{n+2}\right] \times \{2n\}$ and $\bigcup_{k>0} \left(\frac{-1}{n+1}, \frac{-1}{n+1} + \frac{1}{2^{k+1}(n+1)(n+2)}\right] \times \{2n+1+k\}$ are closed in $(-1, 1)$.

Let $n > 1$. Then

$$\left(\frac{-1}{n+1}, \frac{-1}{n+2}\right] \times \{2n\} \cup \bigcup_{k>0} \left(\frac{-1}{n+1}, \frac{-1}{n+1} + \frac{1}{2^{k+1}(n+1)(n+2)}\right] \times \{2n+1+k\} \subseteq A.$$

But

$$\bigcup_{k=0} \left[\frac{-1}{n+1} - \frac{1}{2^{k+1}n(n+1)}, \frac{-1}{n+1}\right] \times \{2n+k\} \subseteq U_{-(n-1)} \subseteq A.$$

Hence the half open intervals from the $n$'th step of construction will "glue" to the closed intervals from the $n - 1$'st step of construction, hence the set $A$ is closed.

Picture below presents the sketch of the set $A$.

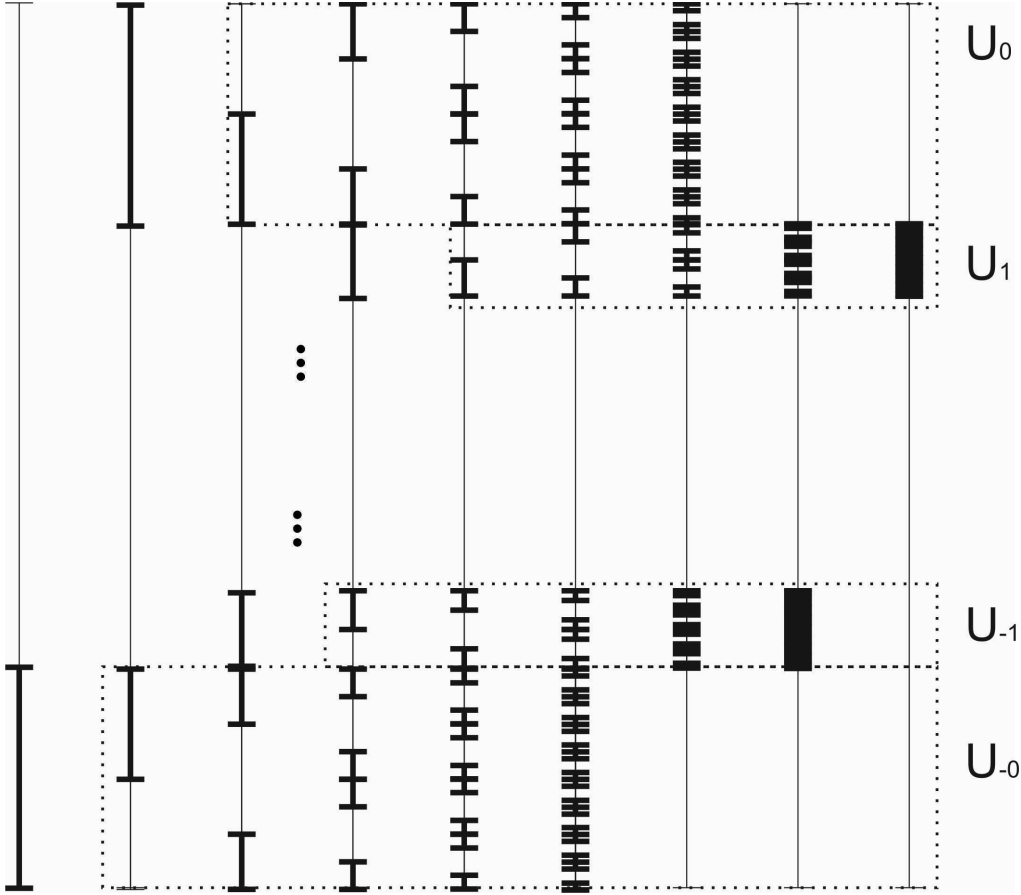

Figure 2:

Set $A$ is closed and as horizontal sections we obtain all sequences except $\left(\overline{x}_n = 0, ..., 0, \underset{n \text{ coordinate}}{1}, 0, 1, 1, 1, 1...\right)$. We will define set $V$ which will be $UU$ by mathematical induction. We do it similarly as in the proof of previous theorem, i.e. we add some point to set $A$ to obtain sequence which we have not covered, but then we loose some other sequence. This lost sequence will be added in some further step.

$n = 0$. Let $V_0 = A \cup \left(\left\{\frac{-1}{3}\right\} x \{0\}\right)$. We obtain sequence $\overline{x}_0 = (1, 0, 1, 1, 1...)$, but we loose sequence $\overline{y}_2 = (0, 0, 1, 1, 1, 1...)$.

$n = 1$. Let $V_1 = V_0 \cup \left(\left\{\frac{-1}{4}\right\} x \{2\}\right) \cup \left(\left\{\frac{-1}{4}\right\} x \{3\}\right)$. We obtain sequence $\overline{y}_2$, but we loose sequence $\overline{y}_4 = (0, 0, 0, 0, 1, 1, 1, 1, 1, ...)$

$n = 2$. Let $V_2 = V_1 \cup \left(\left\{\frac{1}{3}\right\} x \{1\}\right)$. We obtain sequence $\overline{x}_1 = (0, 1, 0, 1, 1, 1, ...)$, but we loose sequence $\overline{y}_3 = (0, 0, 0, 1, 1, 1, 1, ...)$

And we proceed. In each step we obtain some sequence and we loose some sequence. The lost sequences are obtained in further steps, so, in the end, as horizontal sections of the set $V = \bigcup_{n \in \omega} V_n$ we will obtain all sequences from $2^\omega$. The set $V$ is closed, because the set $A$ was closed, and to each of the vertical sections we will add only finitely many points. From the properties of set $A$ and by the construction of $V$ it follows that each sequence is obtained in a unique way. Hence $V$ is UU closed set for $\mathbb{R} \times \omega$. ∎

**Acknowledgement 5** *The author would like to thank Professors Marek Balcerzak and Szymon Głąb for the helpful comments which improved the quality of the paper.*

ALICJA KRZESZOWIEC
E-MAIL: ALICJAKRZESZOWIEC@GMAIL.COM
INSTITUTE OF MATHEMATICS
ŁÓDŹ UNIVERSITY OF TECHNOLOGY
UL. WÓLCZAŃSKA 215
90-924 ŁÓDŹ
POLAND